\documentclass[12pt,reqno]{amsart}

\usepackage{latexsym}

\newcommand{\E}{{\mathbf E}}

\newcommand{\szego}{Szeg\"o }
\newcommand{\nhat}{\raisebox{2pt}{$\wh{\ }$}}

\newcommand{\kahler}{K\"ahler }
\newcommand{\sqrtn}{\sqrt{N}}
\newcommand{\wt}{\widetilde}
\newcommand{\wh}{\widehat}
\newcommand{\PP}{{\mathbb P}}
\newcommand{\R}{{\mathbb R}}
\newcommand{\C}{{\mathbb C}}

\newcommand{\CP}{\C\PP}
\renewcommand{\d}{\partial}
\newcommand{\dbar}{\bar\partial}

\renewcommand{\H}{{\mathbf H}}
\newcommand{\half}{{\frac{1}{2}}}
\newcommand{\vol}{{\operatorname{Vol}}}

\renewcommand{\phi}{\varphi}

\newcommand{\ccal}{\mathcal{C}}

\newcommand{\hcal}{\mathcal{H}}
\newcommand{\lcal}{\mathcal{L}}
\newcommand{\mcal}{\mathcal{M}}

\newcommand{\rcal}{\mathcal{R}}

\newcommand{\al}{\alpha}

\newcommand{\ga}{\gamma}

\newcommand{\la}{\lambda}
\newcommand{\ep}{\varepsilon}
\newcommand{\de}{\delta}

\newcommand{\om}{\omega}

\newcommand{\di}{\displaystyle}
\newtheorem{theo}{{\sc Theorem}}[section]
\newtheorem{cor}[theo]{{\sc Corollary}}
\newtheorem{lem}[theo]{{\sc Lemma}}

\newenvironment{rem}{\medskip\noindent{\it Remark:\/} }{\medskip}

\title[Random polynomials and  Levy concentration of measure]
{Random polynomials of high degree and Levy concentration of
measure}

\author{Bernard Shiffman}
\author{Steve Zelditch}
\address{Department of Mathematics, Johns Hopkins University, Baltimore,
MD
21218, USA}
\email{shiffman@math.jhu.edu, zelditch@math.jhu.edu}

\thanks{Research partially supported by NSF grants DMS-0100474 (first
author) and  DMS-0071358 (second author).}

\date{March 24, 2003}

\begin{document}

\begin{abstract} We show that the ${\mathcal L}^p$ norms  of random sequences of
holomorphic sections $s_N\in H^0(M,L^N)$ of powers of a positive line bundle $L$ over a
compact K\"ahler manifold $M$
 satisfy $$\|s_N\|_p/\|s_N\|_2 =\left\{\begin{array}{ll}
O(1)& \mbox{ for }\
2\le p <\infty\\ O(\sqrt{\log N}) &  \mbox{ for }\  p=\infty  \end{array}\right\} \ \
\mbox{almost surely}.$$ This estimate also holds for  almost-holomorphic sections
of positive line bundles on symplectic manifolds (in the sense of our previous work)
and we give almost sure bounds for the ${\mathcal C}^k$ norms.  Our methods
involve asymptotics of Bergman-Szeg\"o kernels and the concentration of measure
phenomenon.
\end{abstract}

\maketitle

\tableofcontents

\section{Introduction}  Levy concentration of measure occurs when Lipschitz  continuous
functions $f$ on a metric probability  space $(X, d, \mu)$ of large dimension $d$
 are highly concentrated around their median values $\mcal_f$.  In the
fundamental  case
 where $X$ is the unit $N$-sphere $S^N$ with the usual
distance function, and $\mu$ is the
SO$(N +1)$-invariant probability measure,
  the concentration of measure inequality says that
\begin{equation} Prob \left\{ x \in S^{N} : |f(x) - \mcal_f| \geq r \right\} \leq
\exp \left(-\frac  {(N - 1) r^2 }{2 \|f\|_{Lip}^2}\right),
\label{Levy}\end{equation}
where $$\|f\|_{Lip} = \sup_{d(x, y) > 0} \frac{|f(x) - f(y)|}{d(x,
y)|}
$$ is the Lipschitz norm. (See, e.g.
\cite{L}.)

The purpose of this paper is to apply the concentration of measure
inequality for high-dimensional spherical and (associated)
Gaussian ensembles  to determine the distribution of $\lcal^p$
norms of random complex polynomials and more general holomorphic
sections of positive line bundles over complex manifolds. In each
case, we have a sequence of (finite dimensional) Hilbert spaces
${\mathcal H}_N$ of dimensions $d_N \to \infty.$ We will use a
natural inner product $\langle, \rangle_N$ on ${\mathcal H}_N$ to
define
 the unit sphere $S{\mathcal H}_N\subset{\mathcal
H}_N$ endowed with Haar probability measure $\nu_N$.  We also
consider the closely related Gaussian measure $\gamma_N$ on
$({\mathcal H}_N, \langle, \rangle_N)$. In our applications,
${\mathcal H}_N$ will belong to one of the following classes:
\medskip

\noindent{\it I. Complex Ensembles}

\begin{enumerate}

\item [i)] The space ${\mathcal P}_{\C N}^m$ of complex holomorphic
polynomials of degree $N$, with the usual Fubini-Study inner
product;

\item [ii)] More generally, the spaces   $H^0(M, L^N)$
of holomorphic sections of high  powers of a positive hermitian
line bundle $L$, with the inner product
induced by the Hermitian metric;

\item [iii)] Even more generally, the spaces $H^0_J(M, L^N)$ of
almost-holomorphic sections of an ample line bundle $L \to M$ over
a symplectic almost complex manifold.

\end{enumerate}
\medskip

The generalization to almost complex symplectic manifolds is
motivated by the role that almost holomorphic sections have played
in symplectic geometry since the paper of Donaldson \cite{DON}.
Almost holomorphic sections behave  less `deterministically' than
holomorphic sections do, e.g. their zeros sets may or may not be
symplectic submanifolds. This suggests studying them
probabilistically.  In \cite{SZ2} we developed the analytic tools
sufficient for this purpose,  and  initiated a probabilistic study
in \cite{SZ3}. The results of this paper on the almost complex
case are relevant to, and give a natural continuation of, these
earlier articles.

 With certain modifications, our methods also apply
to:
\medskip

\noindent{\it II. Real Ensembles}

\begin{enumerate}

\item  [i)] The space $(V_N(S^m)$ of spherical harmonics of  degree
$N$, with the standard inner product;

\item [ii)] The space ${\mathcal P}_{\R N}^m$ of real polynomials of
degree $N$, with the `Fock space' inner product;

\item [iii)] The space ${\mathcal E}_{\lambda, M, g}$ of linear
combinations of eigenfunctions of the Laplacian on a compact
Riemannian manifold $(M, g)$ with eigenvalue of $\sqrt{\Delta}$ in
the interval $[\lambda, \lambda + 1].$

\end{enumerate}

The pattern  of results and techniques is similar  in the real and
complex cases.  The main difference is in the properties of their
reproducing (Szeg\"o) kernels. The theory of   random real
ensembles is not as developed as  the complex case, and  we plan
to devote a separate article for that purpose. Some of the
relevant steps have already been taken in  \cite{N, V, Ze1}.

The main functionals we consider are the   norms on $S\hcal_N$:
$${\mathcal L}^p (s) = \| s\|_{p} \quad (2\le p\le\infty)\;,\qquad {\mathcal
L}^\infty_{(k)} (s) =
\|\nabla^k s\|_\infty \quad (k=1,2,3,\dots)\;.$$
 We separate
out the cases $p = \infty, p < \infty$ since the
proofs are somewhat different. We also separate out the cases $ k
= 0, k
> 0$ since the case $ k = 0$ is more elementary.

Our first result gives concentration inequalities for sup norms.
The same results hold for both real and complex ensembles, but we
only carry out the proof in the complex case.  We write elements
of ${\mathcal H}_N$ as $s_N$.

\begin{theo}\label{Largesupnorm}  For each of the above complex ensembles,
there exist constants $C > 0$ such that:

 $$  \nu_N\left\{s_N\in S{\mathcal H}_N: \sup_M|s_N|>C\sqrt{\log
N}\right\} < O\left(\frac{1}{N^2}\right)\,.  $$

In fact, for any $k>0$, we can bound the
probabilities by $O(N^{-k})$ by choosing $C$ to be sufficiently
large.

\end{theo}

As a corollary we obtain  almost sure
 bounds on the growth of  $\lcal^\infty$  norms  for
 independent random sequences of $\lcal^2$-normalized holomorphic sections.
To state the result, we introduce the probability sequence space
${\mathcal S} = \prod_{N=1}^{\infty} S{\mathcal H}_N$ with the
measure $\nu = \prod_{N=1}^{\infty} \nu_N$.
The estimate of Theorem  \ref{Largesupnorm} immediately implies
that
$$\limsup_{N\to\infty}\frac{\sup_X|s_N|}{\sqrt{\log N}} \le C\qquad
\mbox{\rm almost surely}\;.$$ Hence we have:

\begin{cor} \label{CX} Sequences of sections  $s_N \in S{\mathcal H}_N$ satisfy:
$$\|s_N\|_{\infty}= O(\sqrt{\log
N})\;\;\; \mbox{ almost surely}.$$

\end{cor}

Results of this type on sup norms were first proved by
Salem-Zygmund \cite{SZ} in the case of random trigonometric
polynomials on the circle , and by Kahane \cite{K} for random
trigonometric polynomials on tori.  Vanderkam \cite{V} generalized
the results to   the case of random spherical harmonics by a
geometric method
 that seems special to the sphere. Nonnenmacher-Voros \cite{NV} obtained
bounds on sup norms for random theta functions on elliptic curves
using properties of the \szego kernel in that setting. Neuheisel
\cite{N} adapted their method to simplify the sup norm estimates
of \cite{V} on random spherical harmonics.  The contribution of
this paper is to give simple and  general results by using
properties of \szego kernels and methods related to Levy
concentration of measure.

Our second result gives  $\lcal^p$ bounds on such sequences.

\begin{theo}\label{largepnorm}  Let $\dim M = m$ and let  $2\le p < \infty$.  Then for
each complex  ensemble, the median values of the ${\mathcal L}^p$
norm on the unit spheres $S\hcal_N$ are bounded by a constant $\al=\al(p,m)$, and
$$  \nu_N \{s_N\in S{\mathcal H}_N: {\mathcal L}^p(s_N) >  r +\al \} \leq
\exp( - Cr^2 N^{ 2 m /p})\;,
$$ for some constant $C>0$. Hence,
sequences of sections $s_N \in S{\mathcal H}_N$ satisfy $\|s_N\|_p
= O(1)$ almost surely.
\end{theo}

Our final results pertain to $\ccal^k$ norms. Almost sure
estimates on $\ccal^k$ norms must take into account the
off-diagonal behavior of the reproducing kernels as well as the
on-diagonal behavior. One of our motivations here is to show that
sequences of increasing degree of almost holomorphic sections of
ample line bundles over symplectic manifolds almost always have
properties similar to the asymptotically holomorphic sections of
Donaldson \cite{DON}. One of these latter properties is a nearly
bounded $\lcal^\infty$ norm. The following result shows that
almost surely, the sup norms of $\lcal^2$-normalized almost
holomorphic sections satisfy comparable estimates to
asymptotically holomorphic sections.

\begin{theo} \label{CK}  For all the ensembles, we have:
$$  \nu_N\left\{s_N\in S{\mathcal H}_N: \sup_M|\nabla^k\; s_N|>C\sqrt{N^k \log
N}\right\} \leq k_N N^{-C^2 N^{-m} (d_N - 1)}.  $$ The
probabilities are thus bounded by $O(N^{-k})$ for any $k$ by
choosing $C$ large enough. Thus sequences $s_N\in S\hcal_N$ satisfy \begin{itemize}
\item $\|\nabla^k
s_N\|_{\infty} = O(\sqrt{N^k \log N})$ almost surely.
\end{itemize}
If $(M, L, \omega)$ is symplectic, then sequences $s_N \in SH^0_J(M,
L^N)$ of almost-holomorphic sections   additionally satisfy:
\begin{itemize}

\item  $\|\bar{\partial}
s_N\|_{\infty} = O(\sqrt{\log N})$ almost surely.

\item  $\|\nabla^k\bar{\partial}
s_N\|_{\infty} = O(\sqrt{N^k\log N})$  almost surely.

\end{itemize}

\end{theo}

 We can also endow $({\mathcal H}_N, \langle,
\rangle_N)$ with
the Gaussian probability measure $\ga_N$ given by: \begin{equation}
\label{ONS} d\gamma_N(s) =  \frac 1{\pi^{d_N}}\,\exp\left(- \sum_{j=1}^d|c_j|^2\right)
\,dc,\qquad
 s = \sum_{j =
1}^{d_N} c_j S_j^N\;, \end{equation}
Where $\{S_j^N\}$ is an orthonormal basis of $\hcal_N$.
Thus, the coefficients $c_j$ of (\ref{ONS}) are independent
complex Gaussian random variables satisfying:
$$\E(c_j) = 0,\quad  \E(c_j c_k) =0,\quad \E(c_j \bar{c}_k) = \delta_{jk},$$
where $\E$ denotes the expected value. It follows immediately from our above results
that for sequences $\{s_N\}\in \prod_{N=1}^\infty \hcal_N$,  with probability measure
$\ga=\prod_{N=1}^\infty \ga_N$, we have:

\begin{itemize}
\item $\|s_N\|_{\infty}/\|s_N\|_2 = O(\sqrt{\log
N})$\quad almost surely.

\item $\|s_N\|_p/\|s_N\|_2 = O(1)$\quad almost surely, for $2\le p<\infty$.

\item $\|\nabla^k
s_N\|_{\infty}/\|s_N\|_2 = O(\sqrt{N^k \log N})$\quad almost surely.
\end{itemize}
\begin{itemize}

\item  $\|\bar{\partial}
s_N\|_{\infty}/\|s_N\|_2 = O(\sqrt{\log N})$\quad almost surely.

\item  $\|\nabla^k\bar{\partial}
s_N\|_{\infty}/\|s_N\|_2 = O(\sqrt{N^k\log N})$\quad  almost surely.

\end{itemize}
(The last two statements are vacuous in the holomorphic case.)

We close the introduction with some open problems that seem of
interest in this area. First, we have assumed throughout that the
spaces of polynomials and sections are equipped with Gaussian
measures. It would be interesting to know how the results would
change if one used other measures, e.g. measures of the form $e^{-
S(f)} {\mathcal D}f$ where $S(f) = \|\nabla f\|^2_{\lcal^2} +
\|f\|_{\lcal^2}^2 + \beta \|f\|_{\lcal^2}^4$ that  arise in quantum
field theory. These are difficult to analyze as $N \to \infty$
since little is known about the
 minimal value of $S$ or the number of its critical points.
In the case of   $S = {\mathcal L}^{\infty}$, it is known
\cite{Bourg1, Bourg2} that the minimal value is bounded as $N \to
\infty$ in the case of $SH^0(\CP^1, {\mathcal O}(N)).$ Bourgain
used Rudin-Shapiro sequences and an estimate on the $\lcal^p$
mapping norm of the \szego kernel. We show in Lemma
\ref{BOURGAIN} that the latter estimate holds for any line
bundle over any \kahler manifold.

\section{Notation and background}

The study of sup-norms of random polynomials has a long history.
Among the earliest articles is that of Paley-Wiener-Zygmund
\cite{PWZ}.  In \cite{SZ}, Salem-Zygmund studied sup-norms of
random trigonometrical polynomials
$$P(t) = \sum_0^N \epsilon_n a_n \cos(n t + \phi).$$
Their main result is a prototype for the subsequent results:
$$Prob\left\{\|P\|_{\infty} <\textstyle \lambda \left[\sum_{n = 1}^N a_n^2\right ]
\sqrt{\log N}\right\} \to 1, \quad N \to \infty.$$ Such estimates
are further developed in the book of Kahane \cite{K}, where
further references may be found. In particular, Kahane (\cite{K},
Section 3.2, Theorem 3) proves a many-variable generalization of
Salem-Zygmund's theorem that is a model for our results: Let
$$P(t_1, \dots, t_m) = \sum  \xi_n f_n(t_1, \dots, t_m)$$
where $\{f_n\}$ are complex trigonometric polynomials in $m$
variables of degrees $\leq N$, and where $\xi_n$ are normal random
variables. Then:
$$Prob\left\{\|P\|_{\infty} \geq\textstyle C s\left [\sum_{n = 1}^N |a_n|^2\right]
\sqrt{\log N} \right\} \leq N^{-2} e^{-s}. $$

Our estimates on sup norms follow however a different path, and
involve scaling asymptotics of Bergman-\szego kernels and
concentration of measure estimates, applied to the spaces $\hcal_N$.  To describe these
Bergman-\szego kernels, we regard the spaces $\hcal_N$ as subspaces of  the
Hilbert space $\lcal^2(X)$ of complex-valued square-integrable
functions on a manifold $X$  ($=$ the $S^1$ bundle
associated to the line bundle $L$). We
denote by
$$\Pi_N: \lcal^2(X) \to \hcal_N$$ the orthogonal
projection onto the subspace; it is given by the Bergman-\szego kernel
$$\Pi_N(x, y) = \sum_{j = 1}^{d_N} S^N_j(x)
\overline{S_j^N(y)}\;,\;\;\;\;\;\; d_N = \dim \hcal_N, $$ where $\{S_j^N\}$ is an
orthonormal basis of $\hcal_N$ (see \S \ref{s-kernels}).

The following well-known elementary probability lemma is central to our
arguments:

\begin{lem}\label{calculus} Let $A\in S^{2d-1}\subset\C^d$, and give
$S^{2d-1}$ Haar probability measure.  Then
$$Prob\left\{P\in S^{2d-1}: |\langle
P,A\rangle|>\la\right\} \le e^{-(d-1)\la^2}\;.$$
\end{lem}

\begin{proof} We can assume without loss of generality that
$A=(1,0,\ldots,0)$.  Let
$$V_\la=\vol \big(\{P\in
S^{2d-1}:|P_1|>\la\}\big)\qquad (0\le\la<1)\,,$$ where $\vol$
denotes $(2d-1)$-dimensional Euclidean volume.  Our desired
probability equals $V_\la/V_0$. Let $\sigma_n=\vol
(S^{2n-1})=\frac{2\pi^n}{(n-1)!}$. We compute
\begin{eqnarray*}V_\la
&=& \int_\la^1\sigma_{d-1}(1-r^2)^{\frac{2d-3}{2}}\frac{2\pi
rdr}{\sqrt{1-r^2}}\ =\
2\pi\sigma_{d-1}\int_\la^1(1-r^2)^{d-2}rdr\\&=&
\frac{\pi\sigma_{d-1}}{d-1}(1-\la^2)^{d-1}\ =\ \sigma_d
(1-\la^2)^{d-1}\,.\end{eqnarray*} Therefore, $$Prob\left\{P\in
S^{2d-1}: |\langle
P,A\rangle|>\la\right\}=V_\la/V_0=(1-\la^2)^{d-1} \le
e^{-(d-1)\la^2}\;.$$ \end{proof}

\begin{rem} In the real case,  from the Levy
concentration of measure inequality (\ref{Levy}) with $f(X)= \Re\langle
P,A\rangle$, we obtain the analogous result:
$$Prob\left\{P\in S^{2d-1}: |\Re\langle
P,A\rangle|>\la\right\} \le e^{-(d-1)\la^2}\;.$$
In the complex case, we then have \begin{eqnarray*}&&Prob\left\{P\in S^{2d-1}:
|\langle P,A\rangle|>\la\right\} \\&&\qquad \le Prob\left\{P\in S^{2d-1}: |\Re\langle
P,A\rangle| >\la/\sqrt 2\right\}+Prob\left\{P\in S^{2d-1}: |\Im\langle
P,A\rangle|>\la/\sqrt 2\right\}\\
&&\qquad\le 2e^{-(d-1)\la^2/2}\;,\end{eqnarray*}
which is not as sharp as Lemma \ref{calculus}.\end{rem}

\subsection{Complex ensembles}

Let us now specify the set-up on \kahler manifolds. We only give a
brief discussion and refer the reader to \cite{BSZ1, BSZ2, SZ1,
SZ2, Ze} for further details.

The simplest example is where ${\mathcal H}_N = {\mathcal
P}_{\C}^N$, the space of holomorphic polynomials
$$f(z_1, \dots, z_m) = \sum_{\alpha \in {\bf N}^m: |\alpha| \leq N}
c_{\alpha} z^{\alpha}.$$ Such polynomials may be regarded as
sections $s \in H^0(\CP^m, {\mathcal O}(N))$ of the $N$-th power
of the hyperplane section bundle ${\mathcal O}(1) \to \CP^m$. The
identification is established by homogenizing each monomial
$$z_1^{\alpha_1} \cdots z_m^{\alpha_m}  \to z_0^{N - |\alpha|} z_1^{\alpha_1} \cdots
z_m^{\alpha_m}$$ to be homogeneous of degree $N$. The
homogenization $\hat{f}$ of $f$ is then determined by its
restriction to $S^{2m + 1} \subset \C^m$ where it is equivariant
under the natural $S^1$ action of $S^{2m + 1} \to \CP^m$:
$\hat{f}(e^{i \theta } z) = e^{i N \theta} \hat{f}(z). $ The space
of such homogenized polynomials will be denoted by ${\mathcal
H}_N$. They satisfy the boundary  Cauchy-Riemann equations
$\dbar_b \hat{f} = 0$ as  boundary values of holomorphic
functions on the unit ball $B_1 \subset \C^m.$

Essentially the same construction exists on any compact algebraic
manifold $(M, \omega)$, i.e. a \kahler manifold  such that
$[\frac{1}{\pi}\omega]$ is an integral cohomology class. There
exists  a hermitian line bundle $(L, h) \to M$ and a metric
connection $\nabla$ on $L$ with curvature given by $\frac{i}{2}
\Theta_L = \omega$. We denote the space of holomorphic sections of
the $N$-th power of $L$ by $H^0(M, L^N)$.

In order to simultaneously analyze sections $s \in  H^0(M, L^N)$
for all $N$,  we lift them to the associated  $S^1$ bundle $$X =
\{v \in L^*: \|v\|_{h^*} =1\} \to M$$ where  $\pi: L^* \to M$
denotes the dual line bundle to $L$ with dual metric $h^*$. We let
$\alpha$ be the connection 1-form on $X$ given by $\nabla$; we
then have $d\alpha =\pi^* \omega$, and thus $\al$ is a contact
form on $X$, i.e., $\al\wedge (d\al)^m$ is a volume form on $X$.

We let $r_{\theta}x =e^{i\theta} x$ ($x\in X$) denote the $S^1$
action on $X$ and denote its infinitesimal generator by
$\frac{\partial}{\partial \theta}$. A section $s$ of $L$
determines an equivariant function $\hat{s}$ on $L^*$ by the rule
$\hat{s}(\lambda) = \left(\lambda, s(z) \right)$ ($\lambda \in
L^*_z, z \in M$).  We  restrict $\hat{s}$ to $X$ to an equivariant
function transforming by  $\hat{s}(r_{\theta} x) = e^{i
\theta}\hat{s}(x)$. Similarly, a section $s_N$ of $L^{N}$
determines an equivariant function $\hat{s}_N$ on $X$: put
\begin{equation} \label{sNhat}\hat{s}_N(\lambda) = \left( \lambda^{\otimes
N}, s_N(z) \right)\,,\quad \la\in X_z\,,\end{equation} where
$\lambda^{\otimes N} = \lambda \otimes \cdots\otimes \lambda$;
then $\hat s_N(r_\theta x) = e^{iN\theta} \hat s_N(x)$. We denote
by $\lcal^2_N(X)$ the space of such equivariant functions
transforming by the $N$-th character, and by ${\mathcal H}_N$ the
subspace of CR functions annihilated by the tangential
Cauchy-Riemann operator $\bar{\partial}_b.$

The space ${\mathcal H}_N$ carries the natural inner product
$$\left\langle \hat{s}, \overline{\hat t} \right\rangle = \int_X \hat s\,
\overline{\hat t} \, dV, \;\;\;\; dV = \alpha \wedge (d
\alpha)^{m-1}. $$ We choose an orthonormal basis $\{S_j^N\}$ and
write every element as
$$\hat{s} = \sum_{j = 1}^{d_N} a_j S^N_j. $$

By the Riemann-Roch formula, we have the estimate for the
dimensions $d_N$:
\begin{equation}\label{RR}
d_N=\frac{c_1(L)^m}{m!}N^m+O(N^{m-1})\,.\end{equation}

\subsubsection{Almost complex ensembles}

This is similar to the complex case except that the complex
structure is non-integrable.  We let $(M, \omega, J)$ be a
compact, almost complex symplectic manifold such that
$[\frac{1}{\pi}\omega]$ is an integral cohomology class, and
choose a hermitian line bundle $(L, h) \to M$ and a  metric
connection $\nabla$ on $L$ with $\frac{i}{2} \Theta_L = \omega$.

  In the general
almost-complex symplectic case it is an almost CR manifold.  The
{\it almost CR structure\/} is defined as follows:
 The kernel of $\alpha$ defines a horizontal hyperplane bundle $H \subset
TX$. Using the projection $\pi: X \to M$, we may lift the
splitting $TM=T^{1,0}M \oplus T^{0,1}M$ to a splitting $H=H^{1,0}
\oplus H^{0,1}$. The almost CR structure on $X$ is defined to be
the splitting
   $TX = H^{1,0} \oplus H^{0,1} \oplus \C \frac{\partial}{\partial
\theta}$.   We  also consider a local orthonormal frame $Z_1,
\dots, Z_n$ of $H^{1,0}$ , resp.\ $\bar{Z}_1, \dots, \bar{Z}_m$ of
$H^{0,1}$, and dual orthonormal coframes $\vartheta_1, \dots,
\vartheta_m,$ resp. $\bar{\vartheta}_1, \dots, \bar{\vartheta}_m$.
On the manifold $X$ we have $d= \d_b +\dbar_b
+\frac{\partial}{\partial \theta}\otimes \alpha$, where
$\partial_b  = \sum_{j = 1}^m {\vartheta}_j \otimes{Z}_j$ and
 $\dbar_b  = \sum_{j = 1}^m \bar{\vartheta}_j \otimes \bar{Z}_j$.
We  define the almost-CR $\bar{\partial}_b$ operator by
$\bar{\partial}_b = df|_{H^{1,0}}$. Note that for an $\lcal^2$
section $s_N$ of $L^N$, we have
\begin{equation}\label{dhorizontal}
(\nabla_{L^N}s_N)\nhat = d^h\hat s_N\,,\end{equation} where
$d^h=\d_b+\dbar_b$ is the horizontal derivative on $X$.

As discussed in \cite{BG, SZ2}, there exists a pseudodifferential
perturbation of $\dbar_b$ which has the main properties of $\dbar_b$ in
the integrable complex case. We denote its kernel by ${\mathcal
H}_N$ and refer to \cite{BG, SZ2} for the definition. By the
Riemann-Roch formula of Boutet de Monvel - Guillemin \cite[\S
14]{BG}, its dimension $d_N$ is again the one in the complex case:
\begin{equation}\label{RRAC}
d_N=\frac{c_1(L)^m}{m!}N^m+O(N^{m-1})\,.\end{equation}

\subsection{Bergman-\szego kernels}\label{s-kernels}

We let $\Pi_N : \lcal^2(X) \rightarrow \hcal_N(X)$ denote the
orthogonal projection.  The Bergman-Szeg\"o kernel $\Pi_N(x,y)$ is
characterized by
\begin{equation} \Pi_N F(x) = \int_X \Pi_N(x,y) F(y) dV_X (y)\,,
\quad F\in\lcal^2(X)\,.
\end{equation} It can be given as
\begin{equation}\label{szego}\Pi_N(x,y)=\sum_{j=1}^{d_N}
S_j^N(x)\overline{ S_j^N(y)}\,,\end{equation} where
$S_1^N,\dots,S_{d_N}^N$ form an orthonormal basis of
$\hcal^2_N(X)$.

The Bergman-\szego kernels determine Kodaira maps
 $\Phi_N : M \to PH^0(M,L^N)'$  to projective
space, defined by
$\Phi_N(z) = \{s_N: s_N(z) = 0\}$. Equivalently, we can choose an
orthonormal basis $S^N_1,\dots,S^N_{d_N}$ of $H^0(M,L^N)$ and
write
\begin{equation}\label{Kmap} \Phi_N : M \to\CP^{d_N-1}\,,\qquad
\Phi_N(z)=\big(S^N_1(z):\dots:S^N_{d_N}(z)\big)\,.\end{equation}
We also define the lifts of the Kodaira maps:
\begin{equation}\label{lift}\wt{\Phi}_N : X \to
\C^{d_N}\,,\qquad
\wt\Phi_N(x)=(S^N_1(x),\dots,S^N_{d_N}(x))\,.\end{equation}

Note that
\begin{equation}\label{PiPhi} \Pi_N(x,y)=\wt\Phi_N(x) \cdot
\overline{\wt\Phi_N(y)}\,;\end{equation} in particular,
\begin{equation}\label{PiPhi2}
\Pi_N(x,x)=\|\wt\Phi_N(x)\|^2\,.\end{equation}

We will need several results on the diagonal and off-diagonal
asymptotics of the Bergman-\szego kernels.  It is proved in
\cite{C, Ze} in the holomorphic case and \cite{SZ2} in the
almost-holomorphic case that there exists a complete asymptotic
expansion:
\begin{equation} \label{CXDIAG}  \Pi_N(z,0;z,0)  =  a_0 N^m +
a_1(z) N^{m-1} + a_2(z) N^{m-2} + \dots \end{equation} for certain smooth
coefficients $a_j(z)$ with $a_0 = \pi^{-m}$.
 Hence, the  maps $\Phi_N$  are well-defined for $N\gg 0$. It follows that
\begin{equation}\label{tyza}\|\wt\Phi_N(x)\|=\Pi_N(x,x)^\half
= \pi^{-m/2} N^{m/2} +O(N^{m/2-1})=(\pi^{-m/2}+\ep_N)N^{m/2}
\,,\end{equation} where $\ep_N$ denotes a term satisfying the
uniform estimate
\begin{equation}\label{epN} \sup_{x\in X}|\ep_N(x)| \le
O\left(\frac{1}{N}\right)\,.\end{equation} As a further corollary
one obtains Tian's almost isometry theorem:
 Let $\omega_{FS}$ denote the Fubini-Study form on $\CP^{d_N-1}$.
Then \begin{equation} \label{TIAN} \|\frac{1}{N}  \Phi_N^*(\omega_{FS}) - \omega\|_{\ccal^k} =
O(\frac{1}{N}) \end{equation} for any $k$.

Off-diagonal asymptotics have been obtained in \cite{SZ2} and have been
studied very precisely in \cite{Ch}. The results are as follows:

\begin{enumerate}

\item [a)] Within a  $\frac{C}{\sqrt{N}}$ neighborhood of the diagonal,
the Bergman-\szego kernel is given by
the scaling asymptotics:
$$\begin{array}{l} N^{-m}  \Pi_N(z_0 + u/\sqrt{N}, \theta/N; z_0 + v/\sqrt{N}, 0)
\sim \Pi_1^{\H}
(u,\theta;  v, 0)\left[1 + O(1/\sqrt{N})\right].
\end{array}$$
Here
$$\Pi^\H_1(u,\theta;v,\phi)= \frac{1}{\pi^m} e^{i(\theta-\phi)+i\Im
(u\cdot \bar v)-\half |u-v|^2}\,$$ is the \szego kernel of the
reduced Heisenberg group.

\medskip
\item [b)] Whenever $d(z, w) \leq C/N^{1/3}$, we have:
\begin{equation}\label{neardiag2} |\Pi_N(z,w)| \le \left(\frac 1{\pi^m}
+o(1)\right) {N^m}\exp\left(-\frac {1-\ep} 2 N
d(z,w)^2\right)+O(N^{-\infty})\;. \end{equation}

\item [c)]  On all of $M$, we have:
\begin{equation}\label{offdiag} |\Pi_N(z, w)| \le CN^m \exp\left(-\la \sqrt N\,
d(z,w) \right)\;.\end{equation}

\end{enumerate}

The near-diagonal scaling asymptotics in (a) is just the first two
terms of a complete asymptotic expansion. Let $P_0\in M$ and
choose a Heisenberg coordinate chart about $P_0$ in the sense of
\cite{SZ2}. Then  \cite[Theorem~3.1]{SZ2}\begin{equation}
\label{neardiag}
\begin{array}{l}
N^{-m}\Pi_N^{P_0}(\frac{u}{\sqrtn},\frac{\theta}{N};
\frac{v}{\sqrtn},\frac{\phi}{N})\\ \\ \qquad =
\Pi^\H_1(u,\theta;v,\phi)\left[1+ \sum_{r = 1}^{K} N^{-r/2}
b_{r}(P_0,u,v) + N^{-(K +1)/2} R_K(P_0,u,v,N)\right]\;,\end{array}
\end{equation}  where  $\|R_K(P_0,u,v,N)\|_{\ccal^j(\{|u|\le \rho,\ |v|\le \rho\}}\le
C_{K,j,\rho}$ for $j\ge 0,\,\rho>0$ and $C_{K,j,\rho}$ is
independent of the point $P_0$ and choice of coordinates.

The estimate (b) on the larger $N^{-1/3}$ balls is  from
\cite[Lemma~5.2(ii)]{SZ2}. The off-diagonal estimate (c) follows by an Agmon distance
argument, as noted by M. Christ
\cite{Ch}; see
\cite[Theorem~2.5]{Be} for an elementary proof.

\section{$\lcal^\infty$ norms: Proof of Theorem \ref{Largesupnorm}}

The proof of   Theorem \ref{Largesupnorm}  is the same in the
complex and almost complex ensembles.

Throughout this section we assume that $\|s_N\|_{\lcal^2} = 1$.
Our aim is to prove:
$$  \nu_N\left\{s_N\in SH^0(M,L^N): \sup_M|s_N|>C\sqrt{\log
N}\right\} < O\left(\frac{1}{N^2}\right)\,, $$ for some constant
$C<+\infty$. (In fact, for any $k>0$, we can bound the
probabilities by $O(N^{-k})$ by choosing $C$ to be sufficiently
large.)

\begin{proof}

Recalling (\ref{lift}), we note that
\begin{equation}\label{coherentstate} \Pi_N(x,y)=\sum_{j=1}^{d_N}
S^N_j(x)\overline{S^N_j(y)}=\langle
\tilde\Phi_N(x),\tilde\Phi_N(y)\rangle\,.
\end{equation}
Let $s_N=\sum_{j=1}^{d_N}c_jS^N_j\ $ ($\sum|c_j|^2=1$) denote a random
element of $SH^0(M,L^N)=S\hcal^2_N(X)$, and write $c=(c_1,\dots,c_{d_N})$.
Recall that
\begin{equation}\label{sNx}s_N(x)=\int_X\Pi_N(x,y)s_N(y)dy=
\sum_{j=1}^{d_N}c_j  S^N_j(x)=
c\cdot\tilde\Phi_N(x)\,.\end{equation} Thus
\begin{equation}\label{cos}|s_N(x)|=\|\tilde\Phi_N(x)\|\cos \theta_x\,,\quad
\mbox{where\ } \cos \theta_x = \frac{\left| c\cdot
\tilde\Phi_N(x)\right|}{\| \tilde\Phi_N(x)\|} \,.\end{equation} (Note that
$\theta_x$ can be interpreted as the distance in $\CP^{d_N-1}$ between
$[\bar c]$ and
$\tilde\Phi_N(x)$.)

Now fix a point $x\in X$.  By Lemma
\ref{calculus}, \begin{eqnarray}\label{prob1}
\nu_N\left\{s_N:\cos\theta_x\ge C N^{-m/2}\sqrt{\log
N}\right\} &\le & \exp\left(-(d_N-1)\frac{C^2\log
N}{N^m}\right)\nonumber \\ &=&\ N^{-C^2 N^{-m}(d_N-1)}\;.\end{eqnarray}

We can cover $M$ by a
collection of $k_N$ balls $B(z^j)$ of radius
\begin{equation}\label{RN} R_N:=\frac{1}{N^{\frac{m+1}{2}}}\end{equation}
 centered at points
$z^1,\dots,z^{k_N}$, where
$$k_N \le O(R^{-2m})\le
O(N^{m(m+1)})\,.$$
By  (\ref{prob1}), we have
\begin{equation} \label{centers} \nu_N \left\{s_N\in
SH^0_J(M,L^N):\max_{j} \cos\theta_{x^j} \ge C N^{-m/2}\sqrt{\log
N}
\right\}\le k_N N^{-C^2 N^{-m}(d_N-1)}\;,\end{equation}
where $x^j$ denotes a point in $X$ lying above $z^j$.

We shall show below that equation (\ref{centers}) together with
(\ref{tyza}) and (\ref{cos}) implies  that the desired sup-norm
estimate holds at the centers of the small balls with high
probability. To obtain our desired estimate on all of $M$, we
first need to extend (\ref{centers}) to points within the balls.
To do this, we consider an arbitrary point $w^j\in B(z^j)$, and
choose points $y^j\in X$ lying above the points $w^j$.  We must
estimate the distance, which we denote by $\de_{N}^j$, between
$\Phi_N(z^j)$ and $\Phi_N(w^j)$ in $\CP^{d_N-1}$. Letting $\ga$
denote the geodesic in $M$ from $z^j$ to $w^j$, we conclude by
(\ref{TIAN}) that
\begin{eqnarray}\de_{N}^j &\le & \int_{\Phi_{N*}\ga}\sqrt{\om_{FS}}
\ =\ \int_\ga \sqrt{\Phi_{N}^*\om_{FS}}
\ \le \ \sqrt{N}\int_\ga
(1+\ep_N)\big)\sqrt{\om}\nonumber \\
& \le& (1+\ep_N)N^\half  R_N
\ =\  \frac{1+\ep_N}{N^{m/2}}\,.
\label{delta} \end{eqnarray}

By the triangle inequality
in $\CP^{d_N-1}$, we have
$|\theta_{x^j}-\theta_{y^j}|\le\de_N^j$. Therefore by (\ref{delta}),
\begin{equation}\label{triangle}\cos\theta_{x^j} \ge
\cos\theta_{y^j}-\de_N^j
\ge
\cos\theta_{y^j}- \frac{1+\ep_N}{N^{m/2}}\,.\end{equation}
By (\ref{triangle}),
$$\cos\theta_{y^j} \ge \frac{(C+1)\sqrt{\log N}}{N^{m/2}}\Rightarrow
\cos\theta_{x^j} \ge \frac{(C+1)\sqrt{\log N} -(1+\ep_N)}{N^{m/2}}
\ge  \frac{C\sqrt{\log N}}{N^{m/2}}$$ and thus
$$\begin{array}{l}\left\{s_N\in
SH^0_J(M,L^N):\sup \cos\theta \ge (C+1) N^{-m/2}\sqrt{\log
N}
\right\}\\[8pt] \qquad\qquad\qquad \subset \left\{s_N\in
SH^0_J(M,L^N):\max_{j} \cos\theta_{x^j} \ge C N^{-m/2}\sqrt{\log
N}
\right\}\,.\end{array}$$

Hence
by  (\ref{centers}),
\begin{equation} \label{cos*} \nu_N \left\{s_N\in
SH^0_J(M,L^N):\sup \cos\theta \ge (C+1) N^{-m/2}\sqrt{\log N}
\right\}\le k_N N^{-C^2 N^{-m}(d_N-1)}\;.\end{equation}

It follows from (\ref{RR}), (\ref{tyza}), (\ref{cos}) and
(\ref{cos*}) that
\begin{eqnarray*}\nu_N \left\{s_N\in SH^0_J(M,L^N):\sup_M |s_N|\ge  (C+2)
\sqrt{\log N}
\right\}\hspace{1.25in}\\
\le k_N N^{-C^2 N^{-m}(d_N-1)}\le O\left(N^{m(m+1)-\frac{C^2}
{m!+1}}\right)\;.\end{eqnarray*} Choosing $C=(m+1)\sqrt{m!+1}$, we
obtain the desired estimate.

\end{proof}

\begin{rem} An alternate proof of this estimate, which does not depend on
Tian's theorem, is given by the case $k=0$ of the $\ccal^k$
estimate in \S \ref{s-Ck}.
\end{rem}

\subsection{Relation to Levy concentration}

The estimate in Theorem  \ref{Largesupnorm} is very closely
related to Levy's estimate.  The proof shows that

\begin{enumerate}

\item [(i)] $\lcal^\infty_N$ is Lipschitz continuous with norm
$\frac{ N^{m/2}}{\sqrt{\log N}} \leq \|\lcal^\infty_N\|_{Lip} \leq
N^{m/2}$. ;

\item [(ii)] The median of $\lcal^\infty_N$ satisfies:   $\mcal_{\lcal^\infty_N}
\leq C_m  \sqrt{\log N}$ for sufficiently large $N$.

\end{enumerate}

Indeed, Lipschitz continuity follows from  equivalence of norms on
finite dimensional vector spaces. To estimate the Lipschitz norm,
we recall the well-known fact that the $\lcal^2$-normalized
`coherent states' $\Phi_N^w(z) = \frac{\Pi_N(z,
w)}{\sqrt{\Pi_N(w,w)}}$ are the global maxima of $\lcal^\infty_N$
on $SH^0(M, L^N)$, as follows from the Schwartz inequality applied
to the reproducing identity $s(z) = \int_{M} \Pi_N(z,w) {s(w)}
dV(w).$ Moreover, $\|\Phi_N^w(z)\|_{\lcal^\infty} =
\sqrt{\Pi_N(w,w)} \sim N^{m/2}$. It follows that
$$\big| \| s_1 + s_2\|_{\infty} - \|s_1\|_{\infty} \big| \leq 3 N^{m/2}.$$
Now let $s_1 $ have $\lcal^\infty$ norm $\leq C \sqrt{\log N}$ and
let $s_1 = \Phi_N^w$ for some $w$. Then we see that
$$\big| \| s_1 + s_2\|_{\infty} - \|s_1\|_{\infty} \big| \geq
\frac{ N^{m/2}}{\sqrt{\log N}}.$$

 It obviously follows from (i)--(ii)   combined with the  Levy
estimate (\ref{Levy}) that (for any $C > 0$)
\begin{equation} \mu \{ s \in SH^0(M, L^N) : f_N^{\infty}(s) \geq C \sqrt{\log N}  \} \leq
 \exp( - C (d_N - 1) \log N /2 N^m).
\end{equation}
Since $d_N \sim N^m$, this  is essentially the same estimate as in
Theorem \ref{Largesupnorm}.

The question arises to find the true order of magnitude of the
median $\mcal_{\lcal^\infty_N}$. It would seem to be smaller than
$\sqrt{\log N}$.

\section{$\lcal^p$ norms: Proof of Theorem \ref{largepnorm}}

We now consider $\lcal^p$ norms for $p < \infty$. We denote by
${\mathcal L}^p_N: S\hcal_N\to\R^+$ the functional ${\mathcal
L}^p(s_N) = \|s_N\|_{\lcal^p}$, $s_N\in S\hcal_N$. Recall that $\hcal_N=H^0(M, L^N)$ in
the holomorphic case, or more generally, $\hcal_N=H^0_J(M, L^N)$ in the symplectic case.
Theorem
\ref{largepnorm} follows from the Levy concentration of measure inequality (\ref{Levy})
applied to  estimates of the Lipschitz norm of
${\mathcal L}^p_N$ and of its median value, which we give in the following two Lemmas.

\subsection{Estimate of the Lipschitz norm of ${\mathcal L}^p_N$}

The first step in the proof of Theorem \ref{largepnorm} is:

\begin{lem} \label{BOURGAIN} The
Lipschitz norm of ${\mathcal L}^p_N$ in dimension $m$  is  $O( N^{m(1/2 -
1/p)})$.
\end{lem}

\begin{proof}

The main point is to show that
\begin{equation}\label{LP}
\sup_{s_N \in S\hcal_N}
\|s_N\|_{\lcal^p(M)} \le C  N^{m(1/2 -
1/p)}\;,
\end{equation}
and hence $$\left|\|s_N\|_{\lcal^p(M)}- \|\tilde s_N\|_{\lcal^p(M)}\right| \le
\|s_N-\tilde s_N\|_{\lcal^p(M)} \le C  N^{m(1/2 -
1/p)} \|s_N-\tilde s_N\|_{\lcal^2(M)}\;.$$

To prove (\ref{LP}), it suffices to show the following estimate for the $\lcal^p
\to \lcal^q$ mapping norm of $\Pi_N$:

\begin{equation} \label{MAPNORMCX} \|\Pi_N f\|_{\lcal^q(M)} \leq C
 N^{m(1/p - 1/q) }
\|f\|_{\lcal^p(M)}. \end{equation}

We shall apply the Shur-Young inequality which bounds the
norm of an integral operator  $K: \lcal^p \to \lcal^q$  by:
$$\|K\|_{\lcal^p \to \lcal^q} \leq C_p\left[\sup_x \int_{M} |K(x, y)|^r
 d\mu(y)\right]^{1/r},\qquad \frac{1}{r} = 1 - \frac{1}{p} + \frac{1}{q}\;. $$
 We break up the integral
 \begin{eqnarray*}  \int_M |\Pi_N(z, w)|^r dV(w) & = &
\int_{d(z, w) \leq N^{-1/3}}
|\Pi_N(z, w)|^r dV(w)   \\ & +&
\int_{d(z, w) \geq  N^{-1/3}} |\Pi_N(z, w)|^r
dV(w).\end{eqnarray*}

For the first term, we  have by (\ref{neardiag2}),
$$\int_{d(z, w) \leq  N^{-1/3}}
|\Pi_N(z, w)|^r dV(w) \le CN^{mr} \int_{\C^m} e^{-rN|u|^2/4}\,du + O(N^{-\infty})
\le C'N^{mr-m}\;.$$
The second term is
rapidly decaying.  Indeed by (\ref{offdiag}),
$$\int_{d(z, w) \geq  N^{-1/3}}
|\Pi_N(z, w)|^r dV(w) \le O(N^{-\infty})\;,$$
and hence \begin{equation}\label{Lr}\int_M
|\Pi_N(z, w)|^r dV(w)
=O(N^{mr-m})\;.\end{equation}

We then obtain (\ref{MAPNORMCX}) from  the Shur-Young inequality and (\ref{Lr}).
\end{proof}

\begin{rem}
Lemma \ref{BOURGAIN} is sharp.  Indeed, we have

\begin{equation}\label{LPCX}
\sup_{s_N \in \hcal_N}
\frac{\|s_N\|_{\lcal^p(M)}}{\|s_N\|_{\lcal^2(M)}} \sim  N^{m(1/2 -
1/p)}
\end{equation}
To prove the lower bound of (\ref{LPCX}), we let $s_N$ be the coherent state
$$\Phi_N^w(z) :=
\frac{\Pi_N(z,w)}{\|\Pi_N(\cdot, w)\|} = N^{-m/2} \Pi_N(z,w)\;.$$ We
have:
$$\|\Phi_N^w\|_{\lcal^p} \sim N^{-m/2} \left[ \int_{M} |\Pi_N(z,w)|^p\,
dV\right]^{1/p} \sim N^{-m/2} N^{m (1 - 1/p) } = N^{m ( 1/2 - 1/p)}.$$
\end{rem}

\subsection{Estimate of the median}

Unlike the case of the sup norm, we can estimate the median
directly by using Chebychev's inequality.

\begin{lem}\label{mediancx}  Let $\mcal_{\lcal^p_N}$ denote the median of
$\lcal^p_N$. Then there is a constant $\alpha=\al(m,p)$ such that
$$\mcal_{\lcal^p_N} \leq \alpha \quad \forall\ N\ge 1\;.$$
\end{lem}

\begin{proof} By Chebychev, we have
\begin{eqnarray*}  \nu_N \{s_N\in S\hcal_N: \lcal^p_N(s_N) > t \} &
\leq&  \frac{1}{t^p} \E \big((\lcal^p_N)^p\big)  \\
& = & \frac{1}{t^p} \int_{S^{2d_N-1}} \int_X \left|\sum_{j=1}^{d_N}c_j
S_j^N(x)\right|^p dV(x) d\nu_N(c)  \;.
\end{eqnarray*}
Let us write
$$S_j^N(x)= \Pi_N(x,x)^\half\, u_j(x) = \|\tilde\Phi_N(x)\|\,u_j(x) \;,$$ so that
$\sum_{j=1}^{d_N}|u_j(x)|^2\equiv 1$. We then have by (\ref{tyza}),
\begin{eqnarray*}  \int_{S^{2d_N-1}} \int_X \left|\sum_{j=1}^{d_N}c_j
S_j^N(x)\right|^p dV(x) d\nu_N(c)\hspace{-1in}\\& =&\int_X \Pi_N(x,x)^{p/2}
\int_{S^{2d_N-1}}\left|\sum_{j=1}^{d_N}c_j u_j(x)\right|^p\, d\nu_N(c)\, dV(x)\\ &=&
A_{p,d_N}\int_X \Pi_N(x,x)^{p/2}
\,dV(x) \ = \ [C_m{A_{p,d_N}} +o(1)]\, N^{mp/2} \;,
\end{eqnarray*}
where
\begin{equation}
\label{A} A_{p,d}=\int_{S^{2d-1}}|w_1|^p\, d\mu(w) \qquad (\mu(S^{2d-1})=1)\;,
\end{equation} and $C_m$ depends only on $m$.  To compute $A_{p,d}$, we evaluate the
integral:
\begin{eqnarray*} \frac 1 {\pi^d} \int_{\C^d} |z_1|^p e^{-\|z\|^2}\, dz &=& \frac 1 \pi
\int_\C|z_1|^pe^{-|z_1|^2}\,d z_1\ =\ \textstyle\Gamma(\frac p2 +1)\\
&=& \frac{\sigma_d}{\pi^d}\int_0^\infty \int_ {S^{2d-1}} |w_1|^p r^p e^{-r^2} r^{2d-1}\,
d\mu(w)\,dr\\
&=& A_{p,d}\, \frac{\sigma_d}{\pi^d} \,\int_0^\infty  r^{p+2d-1} e^{-r^2}\,dr\\
&=& A_{p,d}\, \frac{\Gamma(d+\frac p2)}{\Gamma(d)}\;.\end{eqnarray*}
Therefore
\begin{equation} A_{p,d} = \textstyle\Gamma(\frac p2 +1)\,\di\frac{\Gamma(d)}
{\Gamma(d+\frac p2)} = [\textstyle\Gamma(\frac p2 +1)+o(1)]\,d^{-p/2}\;.\end{equation}

Recalling that $d_N=\frac {c_1(L)^m}{m!}\, N^m+\cdots$, we then have
\begin{equation}\label{Cheb} \nu_N \{s_N\in S\hcal_N: \lcal^p_N(s_N) > t \}
\le \frac 1 {t^p} C_{m,p}\;,\end{equation} for some constant $C_{m,p}$ depending only
on $m$ and $p$.  Substituting $t=\mcal_{\lcal^p_N}$ into (\ref{Cheb}) so that the left
side equals $\half$, we conclude that the sequence $\mcal_{\lcal^p_N}$ of medians is
bounded.\end{proof}

\section{$\ccal^k$ norms: Proof of Theorem \ref{CK}}\label{s-Ck}

Our first aim is to prove the  estimate
\begin{equation} \|\nabla^k s_N\|_{\infty}/\|s_N\|_2 =
O(\sqrt{N^k \log N}) \;\; \mbox{almost surely}
\label{almostsurelyCkbounded} \end{equation} on complex and almost
complex manifolds. The proof  follows the pattern of the above
sup-norm estimate.

We  pause to summarize and review our notation in \cite{SZ2, SZ3}
for the various differential operators that we use in the complex
case:

\begin{enumerate}
\item[a)] Derivatives on $M$:
\begin{itemize}
\item $\frac{\d}{\d z_j} =
\half \frac{\d}{\d x_j}-\frac{i}{2}\frac{\d}{\d y_j}\,,\quad
\frac{\d}{\d\bar
z_j} = \half \frac{\d}{\d x_j}+\frac{i}{2}\frac{\d}{\d y_j}$;\\[-4pt]
\item $Z_j^M=\frac{\d}{\d z_j} +\sum \bar B_{jk}(z)\frac{\d}{\d\bar z_k},\
\bar Z_j^M=\frac{\d}{\d\bar z_j} +\sum  B_{jk}(z)\frac{\d}{\d
z_k},\ B_{jk}(P_0)=0$,\\ $\{Z_1,\dots,Z_m\}$ is a local frame for
$T^{1,0}M$.
\\[-4pt]\end{itemize}
\item[b)] Derivatives on $X$:
\begin{itemize}
\item $\frac{\d^h}{\d z_j} = \frac{\d}{\d
z_j} -A_j(z)\frac{\d}{\d
\theta} =$ horizontal lift of $\frac{\d}{\d z_j}$, \ $A_j(P_0)=0$;\\[-4pt]
\item $Z_j =$ horizontal lift of $Z_j^M$;\\[-4pt]
\item $d^h=\d_b+ \dbar_b$ = horizontal exterior derivative on $X$.\\[-4pt]
\end{itemize}
\item[c)] Covariant derivatives on $M$:
\begin{itemize}
\item $\nabla:\ccal^\infty(M,L^N\otimes (T^*M)^{\otimes k})\to
\ccal^\infty(M,L^N\otimes (T^*M)^{\otimes (k+1)})$;\\[-4pt]
\item
$\nabla^k=\nabla\circ \cdots\circ \nabla:\ccal^\infty(M,L^N)\to
\ccal^\infty(M,L^N\otimes (T^*M)^{\otimes k} )$;\\[-4pt]
\item $\nabla=\d +\dbar,\ \  \dbar:\ccal^\infty(M,L^N)\to
\ccal^\infty(M,L^N\otimes T^{*0,1}M)$.\\[-4pt]
\end{itemize}
\item[d)] Derivatives on $X\times X$:
\begin{itemize}
\item $d_j^1, d_j^2$: the operator $\frac{\d^h}{\d z_j}$ applied to the first
and second factors, respectively;\\[-4pt]
\item $Z_j^1, Z_j^2$: the operator $Z_j$ applied to the first
and second factors, respectively.
\end{itemize}
\end{enumerate}

\subsection{Derivatives of holomorphic and almost holomorphic sections}

To prove
(\ref{almostsurelyCkbounded}), we first note a consequence (Lemma \ref{horiz-bound2})
of our near-diagonal asymptotics.  Recall that a differential operator on $X$ is
horizontal if it is generated by horizontal vector fields.  In particular the operators
$\nabla^k:\ccal^\infty(M,L^N)\to \ccal^\infty(M,L^N\otimes (T^*M)^{\otimes
k})$ are given by (vector valued) horizontal differential operators
(independent of $N$) on
$X$. By definition, horizontal differential operators on
$X\times X$ are generated by the horizontal differential operators on the
first and second factors. We begin with the following estimate:

\begin{lem}\label{horiz-bound1} Let $P_k$ be a horizontal differential
operator of order $k$ on $X\times X$. Then
$$P_k\Pi_N(x,y)|_{x=y}=O(N^{m+k/2})\,.$$\end{lem}

\begin{proof}
Let $x_0=(P_0,0)$ be an arbitrary point of $X$, and choose local
real `Heisenberg' coordinates $(x_1,\dots,x_{2m},\theta)$ about $(P_0,0)$
as in the hypothesis of Theorem 3.1 of \cite{SZ2} (with $z_q=x_q+ix_{m+q}$). We let
$\frac{\d^h}{\d x_q}$ denote the horizontal lift of $\frac{\d}{\d x_q}$ to $X$:
$$\frac{\d^h}{\d x_q} = \frac{\d}{\d x_q} -\tilde A_q(x)\frac{\d}{\d
\theta}\,,\quad \tilde
A_q =(\al, \frac{\d}{\d x_q})\,.$$ Since $\left.\frac{\d}{\d
x_q}\right|_{x_0}$ is assumed to be horizontal, we have $\tilde A_q(P_0)=0$.

We let $d_q^1, d_q^2$ denote the operator $\frac{\d^h}{\d x_q}$
applied to the first and second factors, respectively, on $X\times
X$. For this result, we need only the zeroth order estimate of
(\ref{neardiag}):
\begin{equation}\label{zeroth}
\Pi_N(\frac{u}{\sqrtn},\frac{s}{N}; \frac{v}{\sqrtn},\frac{t}{N})
= {N^m}e^{i(s-t)+\psi_2(u,v)} \rcal(P_0,u,v,N)\,,\end{equation}
where $\rcal(P_0,u,v,N)$ denotes a term satisfying the remainder
estimate of (\ref{neardiag}):
$$\|\rcal(P_0,u,v,N)\|_{\ccal^j(\{|u|\le \rho,\ |v|\le \rho\}}\le
C_{j,\rho}$$ for $j\ge 0,\,\rho>0$, where $C_{j,\rho}$ is
independent of the point $P_0$ and choice of coordinates. Here,
\begin{equation}  \psi_2(u,v) = u \cdot\bar{v} - \half(|u|^2 +
|v|^2).\end{equation}

Differentiating (\ref{zeroth})  and noting that $\d/\d x_q = \sqrtn
\d/\d u_q$, $\d/\d\theta = N\d/\d s$, we have
\begin{eqnarray}\lefteqn{d_q^1\Pi_N(\frac{u}{\sqrtn},\frac{s}{N};
\frac{v}{\sqrtn},\frac{t}{N})} \nonumber \\
&=& \sqrtn \left( \frac{\d}{\d
u_q}-\sqrtn\tilde A_q(P_0+{\textstyle{\frac{u}{\sqrtn}}})\frac{\d}{\d
s}\right)
\left( {N^m}e^{i(s-t)+\psi_2(u,v)}\rcal\right) \nonumber\\
&=&{N^{m+1/2}}
e^{i(s-t)+\psi_2(u,v)}\left\{\left[L_q(u,v)-i\sqrtn \tilde
A_q({\textstyle{\frac{u}{\sqrtn}}})\right] +\frac{\d}{\d u_q}\rcal
\right\}
\nonumber\\ &=& {N^{m+1/2}}
e^{i(s-t)+\psi_2(u,v)} \wt\rcal \ = \
O(N^{m+1/2})\,,\label{expand0}
\end{eqnarray} where $L_q:=\frac{\d\psi_2}{\d u_q}$ is a linear function.
The same estimate holds for $d^2_q\Pi_N$.  Indeed, the above computation
yields:
\begin{equation}\label{estimate0} d_q^je^{i(s-t)+\psi_2(u,v)}
\rcal(P_0,u,v,N) = \sqrtn
e^{i(s-t)+\psi_2(u,v)} \wt \rcal(P_0,u,v,N)\,,\end{equation} for $j=1,2,\
q=1,\dots,2m$.  The desired estimate follows by iterating
(\ref{estimate0}). \end{proof}

\begin{rem} The assumption that $P_k$ is horizontal in Lemma
\ref{horiz-bound1} is necessary, since the operator $\d/\d\theta$
multiplies the estimate by $N$ instead of $\sqrtn$. \end{rem}

\begin{lem} Let $P_k$ be a horizontal differential operator of
order
$k$ on
$X$. Then
$$ \sup_X\|P_k \tilde\Phi_N\| = O(N^\frac{m+k}{2})\,.$$
\label{horiz-bound2}\end{lem}

\begin{proof} Let $P_k^1,\ P_k^2$  denote the operator $P_k$ applied to
the first and second factors, respectively, on $X\times X$. Differentiating
(\ref{coherentstate}) and restricting to the diagonal, we obtain
\begin{equation}\label{coherentstate1}
P_k^1 \bar P_k^2\Pi_N(x,x)= \left
\|P_k\tilde\Phi_N(x)\right\|^2\,.
\end{equation}  The conclusion follows from (\ref{coherentstate1}) and Lemma
\ref{horiz-bound1} applied to the horizontal differential operator (of order
$2k$)
$P_k^1
\bar P_k^2$ on $X\times X$.
\end{proof}

We are now ready to use the small-ball method of the previous section to
show that $\|\nabla^k s_N\|_{\infty}/\|s_N\|_{\lcal^2} = O(\sqrt{N^k \log
N})$ {\it almost surely.\/}  It is sufficient to show that
\begin{equation} \nu_N\left\{s_N\in SH^0_J(M,L^N): \sup_M|\nabla^k
s_N|>C\sqrt{N^k   \log N}\right\} <
O\left(\frac{1}{N^2}\right)\,,
\label{largeCknorm}\end{equation} for $C$ sufficiently large.
To verify (\ref{largeCknorm}), we may regard $s_N$ as a function on
$X$ and replace
$\nabla^k$ by a horizontal $r_\theta$-invariant differential operator of
order
$k$ on
$X$.

As before, we let $s_N=\sum c_j S^N_j$ denote a random element of
$S\hcal^2_N(X)$.
By (\ref{sNx}), we have
\begin{equation}\label{dqsn}P_ks_N(x)=\int_X
P_k^1\Pi_N(x,y)s_N(y)dy= \sum_{j=1}^{d_N}c_j P_kS^N_j(x) = c\cdot
P_k\tilde\Phi_N(x)
\,.\end{equation} We then have
\begin{equation}\label{cos1}|P_ks_N(x)|=\|P_k\tilde\Phi_N(x)\|\cos
\theta_x\,,\quad
\mbox{where\
} \cos \theta_x = \frac{\left| c\cdot
P_k\tilde\Phi_N(x)\right|}{\|
P_k\tilde\Phi_N(x)\|} \,.\end{equation}
Now fix a point $x\in X$. As before, (\ref{prob1}) holds, and hence by Lemma
\ref{horiz-bound2} we have
\begin{equation} \label{supPk} \nu_N \left\{s_N\in S\hcal^2_N:
|P_ks_N(x)| \ge C'
\sqrt{N^k\log N}
\right\}\le k_NN^{-C^2 N^{-m}(d_N-1)}\;,\end{equation}
where $C'=C \sup_{N,x}N^{-(m+k)/2}|P_k\tilde\Phi_N(x)|$.

We again  cover $M$ by a collection of $k_N$ very small balls $B(z^j)$ of
radius
$R_N=N^{-\frac{m+1}{2}}$ and first show that the probability of the
required condition holding at the centers of all the balls is small.
Choosing points
$x^j\in X$ lying above the centers $z^j$ of the balls, we
then have
\begin{equation} \label{centers2} \nu_N \left\{s_N\in S\hcal^2_N:\max_{j}
|P_ks_N(x_j)| \ge C'
\sqrt{N^k\log N}
\right\}\le k_NN^{-C^2 N^{-m}(d_N-1)}\;.\end{equation}
Now suppose that $w^j$ is an arbitrary point in $B(z^j)$, and let $y^j$ be
the point of $X$ above $w^j$ such that the horizontal lift of the geodesic
from
$z^j$ to $w^j$ connects $x^j$ and $y^j$.  Hence by Lemma
\ref{horiz-bound2}, we have
\begin{equation}\label{dschange}\|P_k\tilde
\Phi_N(x^j)-P_k\tilde
\Phi_N(y^j)\|\le \sup_M\|d^h(P_k\tilde
\Phi_N)\|
r_N = O(N^{\frac{m+k+1}{2}})r_N=O(N^{\frac{k}{2}})\,.\end{equation}
It follows as before from (\ref{centers2}) and
(\ref{dschange}) that
\begin{eqnarray*} \nu_N \left\{s_N\in S\hcal^2_N:\sup_X
|P_ks_N| \ge (C'+1)
\sqrt{N^k\log N}
\right\}\hspace{1.25in}\\
\le k_N N^{-C^2 N^{-m}(d_N-1)} \le
O\left(N^{m(m+1)-\frac{C^2}
{m!+1}}\right)\;.\end{eqnarray*}
(Here, we used the fact that $|P_ks_N|$ is constant on the fibers of
$\pi:X\to M$.) Thus, (\ref{largeCknorm})  holds with $C$ sufficiently
large.\qed

\subsection{$\overline{\partial}$ derivatives of  almost holomorphic sequences}

In this section we obtain additional results on the complex
derivatives of almost holomorphic sections.
The results are of course trivial in the holomorphic case. As mentioned in the
introduction, they are relevant to the use of asymptotically holomorphic sections
in almost complex geometry.

\subsubsection{The estimate $\|\bar{\partial}
s_N\|_{\infty}/\|s_N\|_2 =
O(\sqrt{\log N})$ almost surely}

The proof of the $\dbar s_N$ estimate follows the pattern of the
above estimate.  However, there is one crucial difference: we must
show the following upper bound for the modulus of
$\dbar_b\tilde\Phi_N$. This estimate is a factor of $\sqrtn$
better than the one for $d^h \tilde\Phi_N$ arising from Lemma
\ref{horiz-bound2}; the proof depends on the precise second order
approximation of Theorem 3.1 of \cite{SZ2} (see (\ref{neardiag})).

\begin{lem}\label{tyzabar} $\qquad\sup_X\|\dbar_b\tilde\Phi_N(x)
\|\le O(N^{m/2})
\,.$\end{lem}
\begin{proof}
Let $x_0=(P_0,0)$ be an arbitrary point of $X$, and choose
preferred local coordinates $(z_1,\dots,z_{m},\theta)$ about
$(P_0,0)$ as in the hypothesis of Theorem \ref{neardiag}.  We lift
a local frame $\{\bar Z_q^M\}$  to obtain the local frame $\{\bar
Z_1,\dots,\bar Z_m\}$ for $H^{0,1}X$ given by
\begin{equation}\label{barz}\bar Z_q = \frac{\d^h}{\d\bar z_q} +
\sum_{r=1}^m B_{qr}(z)\frac{\d^h}{\d z_r}\,,\qquad
B_{qr}(P_0)=0\,.\end{equation}
It suffices to show that
\begin{equation}\label{barzh}N^{-m/2}|\bar Z_q \tilde \Phi_N (x_0)| \le
C\,,\end{equation} where $C$ is a constant independent of $x_0$.

By (\ref{neardiag}), we have
\begin{eqnarray} \lefteqn{N^{-m}
\Pi_N(\frac{u}{\sqrtn},\frac{s}{N};
\frac{v}{\sqrtn},\frac{t}{N})}\nonumber \\& =&
\frac{1}{\pi^m}\phi_0(u,s)\overline{\phi_0(v,t)}e^{u\cdot \bar v}
\left[1 + \frac{1}{\sqrt{N}}
b_1(P_0,u,v)+\frac{1}{N}R_2(P_0,u,v,N)\right]
\,,\label{asymptotics1}
\end{eqnarray} where
$$\phi_0(z,\theta)=e^{i\theta-|z|^2/2}\,.$$
(The function $\phi_0$ is the `ground state' for the `annihilation
operators' $\bar Z_q$ in the Heisenberg model; see the remark in
\cite[\S 1.3.2]{BSZ2}).  In our case, $\bar Z_q\phi_0$ does not
vanish as in the model case, but instead  satisfies the asymptotic
bound (\ref{groundstate}) below.)

We have (see \cite{SZ2}),
\begin{equation}\frac{\d^h}{\d\bar z_q} = \frac{\d}{\d\bar z_q}
+\left[ -\frac{i}{2}z_q-R_1^{\bar A_q}(z)\right]\frac{\d}{\d \theta}
\,,\label{dh}\end{equation} where
$R_1^{\bar A_q}(z)=O(|z|^2)$.  Recalling that $z=u/\sqrtn$,
$\theta=s/N$, we note that $\phi_0(u,s)=e^{iN\theta
-N|z|^2/2}=\phi_0(z,\theta)^N$,
and thus by (\ref{dh}),
\begin{equation}\label{groundstate} \frac{\d^h}{\d\bar z_q}
\phi_0 (u,s)=\frac{\d^h}{\d\bar z_q}e^{iN\theta -N|z|^2/2}  = -iN
R_1^{\bar A_q}(\frac{u}{\sqrtn}) \phi_0(u,s) =
\rcal(P_0,u,N)\phi_0(u,s)\,,\end{equation} where as before $\rcal$
denotes a term satisfying the remainder estimate of
(\ref{neardiag}).

We let
$ Z_q^1,  Z_q^2$ denote the operator
$ Z_q$ applied to the first and second factors, respectively,
on $X\times X$; we similarly let $d_q^1,
d_q^2$ denote the operator $\frac{\d^h}{\d z_q}$ applied to the factors of
$X\times X$. Equation (\ref{coherentstate1}) tells us that
\begin{equation}\label{barzh2}
\|\bar Z_q\tilde\Phi_N(x)\|^2=
\bar Z_q^1 Z_q^2 \Pi_N(x,x)\,.
\end{equation}
By (\ref{barz}), \begin{equation}\label{4terms}\bar Z_q^1 Z_q^2
=  \left(\overline{d^1_q}+\sum_{r=1}^mB_{qr}(z)
d^1_r\right)
\left(d^2_q+\sum_{\rho=1}^m \bar B_{q\rho}(w)
\overline{d^2_\rho}\right)
\,,\end{equation} where we recall that $B_{qr}(P_0)=0$.

Differentiating
(\ref{asymptotics1}), again noting that
$\d/\d z_q =\sqrtn\d/\d u_q$, $\d/\d w_q =\sqrtn\d/\d v_q$ and using
(\ref{groundstate}), we obtain
\begin{eqnarray} \ \lefteqn{N^{-m}\left(\overline{ d^1_q} d^2_q
\Pi_N \right)(\frac{u}{\sqrtn},\frac{s}{N};
\frac{v}{\sqrtn},\frac{t}{N})}\nonumber \\
&\quad =&\frac{1}{\pi^m}
\phi_0(u,s)\overline{\phi_0(v,t)}e^{u\cdot
\bar v} \left[\sqrtn \frac{\d^2}{\d\bar u_q\d v_q}b_1 +\wt\rcal\right]
\,.\label{asymptotics2}\end{eqnarray}
Since $b_1$ has no terms that are quadratic in $(u,\bar u,v,\bar v)$, it
follows from (\ref{4terms})--(\ref{asymptotics2}) that
\begin{equation} N^{-m}\left|\bar Z_q^1 Z_q^2
\Pi_N(P_0,0;P_0,0)\right| =
N^{-m}\left|\overline{ d^1_q} d^2_q
\Pi_N(P_0,0;P_0,0)\right| = \frac{1}{\pi^m}\left|
\wt\rcal(P_0,0,0,N)\right| \le O(1)\,. \label{barzh1}\end{equation}
The desired estimate (\ref{barzh}) now follows immediately from
(\ref{barzh2}) and (\ref{barzh1}). \end{proof}

By covering $M$ with small
balls and repeating the argument of the previous section, using Lemma
\ref{tyzabar}, we conclude that
\begin{equation}
\nu_N\left\{s_N\in SH^0_J(M,L^N):
\sup_M|\dbar s_N|>C\sqrt{\log N}\right\} < O\left(\frac{1}{N^2}\right)\,.
\end{equation}
Thus $\|\bar{\partial}
s_N\|_{\infty}/\|s_N\|_2 =
O(\sqrt{\log N})$ almost surely.

\subsubsection {The estimate $\|\nabla^k\bar{\partial}
s_N\|_{\infty}/\|s_N\|_2 =
O(\sqrt{N^k\log N})$ almost surely}

To obtain this final estimate of Theorem \ref{CK}, it suffices to
verify the probability estimate
\begin{equation} \nu_N\left\{s_N\in SH^0_J(M,L^N): \sup_M|\nabla^k \dbar
s_N|>C\sqrt{N^k   \log N}\right\} <
O\left(\frac{1}{N^2}\right)\,.
\label{largeCknorm2}\end{equation}
Equation (\ref{largeCknorm2}) follows by again repeating the argument of
\S \ref{almostsurelyCkbounded}, using the following lemma.

\begin{lem} Let $P_k$ be a horizontal differential operator of
order
$k$ on
$X$ ($k\ge 0$). Then
$$ \sup_X|P_k \dbar_b\tilde\Phi_N| = O(N^\frac{m+k}{2})\,.$$
\label{horiz-bound3}\end{lem}

\begin{proof} It suffices to show that
\begin{equation}\label{supkq}\sup_U|P_k \bar Z_q^k\tilde\Phi_N| =
O(N^\frac{m+k}{2})\end{equation} for a local frame $\{\bar Z_q\}$ of
$T^{0,1}M$ over
$U$. As before, we have
\begin{equation}\label{coherentstate2}
P_k^1  \bar P_k^2 \bar Z_q^1 Z_q^1 \Pi_N(x,x)= \left
|P_k\bar Z_q^{h}\tilde\Phi_N(x)\right|^2\,.
\end{equation}

We claim that
\begin{equation}\label{step0} N^{-m}\bar Z^1_q Z^2_q
 \Pi_N = \frac{1}{\pi^m}
e^{i(s-t)+\psi_2(u,v)}\left[\sqrtn  \frac{\d^2}{\d\bar u_q\d v_q}b_1
+\rcal(P_0,u,v,N)\right]\,.\end{equation}
To obtain the estimate (\ref{step0}), we recall from (\ref{4terms}) in
the proof of Lemma \ref{tyzabar} that
\begin{equation}\label{4terms*}\bar Z_q^1 Z_q^2 =
\overline{d^1_q}d^2_q+
\sum_{\rho=1}^m \bar B_{q\rho}(w)\overline{d^1_q}
\overline{d^2_\rho} + \sum_{r=1}^mB_{qr}(z) d^1_r d^2_q +
\sum_{r,\rho}B_{qr}(z)\bar B_{q\rho}(w)  d^1_r\overline{d^2_\rho}\,.
\end{equation}
Equation (\ref{asymptotics2}) says that the first term of
$N^{-m}\bar Z^1_q Z^2_q
\Pi_N$
coming from the expansion (\ref{4terms*}) satisfies the estimate of
(\ref{step0}).  To obtain the estimate for the second term, we compute:
\begin{eqnarray}
\lefteqn{N^{-m}\overline{d^2_\rho}
\Pi_N(\frac{u}{\sqrtn},\frac{s}{N};
\frac{v}{\sqrtn},\frac{t}{N})\ =\ \frac{\sqrtn}{\pi^m}
e^{i(s-t)+\psi_2(u,v)}}\nonumber \\
& & \cdot  \left(\left[\frac{\d\psi_2}{\d \bar v_\rho}+i\sqrtn
A_\rho(\frac{v}{\sqrtn})\right] \left[1+ \frac{1}{\sqrtn}b_1+
\frac{1}{N}R_2
\right] +\frac{1}{\sqrtn}\frac{\d\ b_1}{\d\bar v_\rho}
+\frac{1}{N} \frac{\d\ R_2}{\d\bar v_\rho}\right)\nonumber \\
&=& \frac{\sqrtn}{\pi^m}
\phi_0(u,s)\overline{\phi_0(v,t)}e^{u\cdot
\bar v} \left[\frac{\d\psi_2}{\d \bar
v_\rho} +L_\rho(v)+\frac{1}{\sqrtn}\wt\rcal\right]\,,\end{eqnarray}
where $ L_{\rho}$ is a linear function. Since $\d^2\psi_2/\d\bar
u_q \d\bar v_\rho \equiv 0$, it then follows that
\begin{equation} N^{-m}\overline{ d^1_q}\overline{d^2_\rho}
 \Pi_N = \frac{\sqrtn}{\pi^m}
e^{i(s-t)+\psi_2(u,v)} \frac{\d}{\d \bar
u_q}\wt\rcal(P_0,u,v,N)\,.\label{mixedterm}\end{equation}
The estimate (\ref{step0}) for the second term follows from
(\ref{mixedterm}), using the fact that
$B_{q\rho}(\frac{v}{\sqrtn})=\frac{1}{\sqrtn}L_{q\rho}(v)+\cdots$.
The proofs of the estimate for the third and fourth terms are
similar.

The desired estimate (\ref{supkq}) follows as before from
(\ref{coherentstate2}), (\ref{step0}), and (\ref{estimate0}), using the
fact that $ \frac{\d^2}{\d\bar
u_q\d v_q}b_1$ is linear.
\end{proof}

\end{document}